\documentclass[a4paper,10pt]{article}
\usepackage{graphicx, graphics, amssymb, amsfonts, amsmath,amsthm}

\newtheorem{theorem}{Theorem}
\newtheorem{corollary}{Corollary}
\newtheorem{lemma}{Lemma}

\title{Ergodicity and Percolation for Variants of One-dimensional Voter Models}
\author{Y. Mohylevskyy$^{a,}$\footnote{Research supported in part by NSF grants OISE-0730136 and MPS-1007524} , C.M. Newman$^{a,}$\footnotemark[1] , K. Ravishankar$^{b}$ \\
\begin{scriptsize}                                                                        
$^{a}$ Courant institute of Mathematical Sciences, New York University, New York, NY 10012, USA 
\end{scriptsize}
\\
\begin{scriptsize} 
$^{b}$ Department of Mathematics, SUNY College at New Paltz, New Paltz, NY 12561, USA
\end{scriptsize}
}

\begin{document}

\maketitle
\begin{abstract}
We study variants of one-dimensional $q$-color nearest-neighbor voter models in discrete time. In addition to the usual voter model transitions in which a color is chosen from the left or right neighbor of a site there are two types of noisy transitions. One is bulk nucleation where a new random color is chosen. The other is boundary nucleation where a random color is chosen only if the two neighbors have distinct colors. We prove under a variety of conditions on $q$ and the magnitudes of the two noise parameters that the system is ergodic, i.e., there is convergence to a unique invariant distribution. The methods are percolation-based using the graphical representation of the model which consists of coalescing random walks combined with branching (boundary nucleation) and dying (bulk nucleation).  
\end{abstract}

\section{Introduction}
In this paper we consider a class of one-dimensional interacting particle systems with a focus on ergodicity, i.e., whether there is always convergence to a unique invariant distribution. The evolving state of the system is an assignment of colors from $\{1,...,q\}$ to the sites in $\mathbb{Z}$ and the transition probabilities of a given site depend only on the color of its left and right neighbors. We actually study a two-parameter, $\epsilon$ and $\delta$, family of models which includes classic voter models with or without noise and stochastic Ising and Potts models at zero and nonzero temperature. Together the two parameters control the rates of bulk nucleation of random colors and boundary nucleation when the two neighbors disagree in color. Our main results are about ergodicity under various conditions on $\epsilon$, $\delta$, and $q$. The methods we use are percolation-based and rely on the fact that the dual system is a model of coalescing random walks that also have branching and dying. 

We start by recalling the classic discrete time voter model in one dimension \cite{Liggett}. The voters or particles are located at each site of $\mathbb{Z}$ and take one of the $q$ possible opinions (or colors). At each point in time $t$ the particles make decisions on whether to keep or change their color. The decision consists of randomly choosing one of the two neighbors and adopting its color from time $t-1$. The process of choosing a neighbor may be represented by drawing a random arrow to one of the nearest neighbors at time $t-1$ and thus the history or genealogy of the color of a particular site may be traced to time zero by the path of backward arrows coming out of each space-time site. The distribution of the genealogy paths is that of coalescing simple symmetric random walks.

Noise can be introduced in the voter model by having each site make a decision as before with probability $p$ or else choose a color uniformly at random  out of the $q$ possible ones with probability $1-p$. The genealogy of the noisy voter model is represented by the coalescing random walks with ``dying'', where a path dies at a point where a color choice was made uniformly at random. Such points, which we call bulk nucleation points, represent death of a genealogy path or birth of an opinion. If the opinions are +1 or -1, the noisy voter model coincides with the stochastic Ising model with Glauber dynamics. 

The models we consider in this paper allow for the number of colors of a single site to be any positive integer $q\geq2$. The decision that each particle makes at each point in time has the two possibilities for a noisy voter model: one is to choose a color of a randomly chosen neighbor and another is  to generate a uniformly random color. A third possibility is to take the color of the two neighbors if their colors agree or to generate a uniformly random color if they disagree (boundary nucleation). A site where such a decision is made we represent by a double-arrow in the genealogy graph. To determine the probability of these three choices we introduce two parameters, $\delta$ and $\epsilon$. This model, for certain choices of $\delta$ and $\epsilon$, is related to stochastic Potts models with $q$ colors (for the relation with zero temperature continuous time Potts models see \cite{Derrida}, with continuous time nonzero temperature Ising models see \cite{Fontes}, with nonzero temperature Potts models in 
continuous and discrete time see \cite{Newman}). In the space-time coordinates the evolution of the system may be represented by coloring vertices of $\mathbb{Z}^{2}$. Since a color of vertex $(z,t)$ depends on colors of vertices $(z-1,t-1)$ and $(z+1,t-1)$ only, there are actually two independent systems evolving: one lives on the sub-lattice with $z+t$ even and another with $z+t$ odd. Thus, we can restrict our lattice to be the part of $\mathbb{Z}^{2}$ which contains only vertices with $z+t$ even. 

The special case where the number of colors $q\rightarrow\infty$ is also considered. In this case it is convenient to view coloring as partitioning sites into equivalence classes of the same color rather than coloring the sites, and at nucleation points the color (or equivalence class) that is supposed to be assigned uniformly at random is chosen to be completely new and different from all the previously existing ones. 

The main issues that we investigate can be described by the term memory loss: to what extent does the system in the long run depend on the initial conditions. The existence of an invariant distribution is known for all lattice systems in discrete and continuous time and is based on compactness of the state space. See \cite{Toom} for discrete time systems and \cite{Liggett} for continuous time systems. If the invariant distribution is not unique, then memory is not lost. Limiting behavior will depend on the initial distribution as there are at least as many possible limits as there are invariant distributions. Ergodicity means convergence to a unique possible limit for all initial conditions and is the strongest degree of memory loss. We prove ergodicity under certain conditions on the parameters $\delta$ and $\epsilon$. For others the question remains open. For all nonzero parameters we prove a weaker type of memory loss -- color permutation invariance of any limit or sub-sequence limit distribution. This 
shows that if one begins with all sites of the same color, any sub-sequence limit distribution does not depend on the starting color.

There are natural diffusive scaling limit analogues of the model considered in this paper where space-time is continuous rather than discrete and random walks become Brownian motions. Some of these, such as the Brownian web and net and related continuum voter models, have already been studied \cite{Fontes} and others are in preparation \cite{Newman}. One reason we focus here on the discrete time model is because the only existing published result \cite{Sun} about convergence to these limits starts from the discrete time lattice model. Another reason is to allow us to compare ergodicity and percolation properties of the model -- see the Remark following the proof of Theorem 2 in Section 3 below. We expect that the discrete model results of this paper can be extended to the continuum setting, but do not explore that here. We note however that there are a number of papers \cite{Athreya,Cox1,Cox2,Griffeath} that treat related continuous time lattice models that also have a natural percolation substructure to 
their dual models. 

A formal definition of the model and presentation of main results is given in the next section. Section 3 contains all the proofs. The proofs are based on the underlying directed arrow percolation process which traces the genealogy of colors. Percolation methods such as dynamic renormalization (see \cite{Barsky}) or enhancement (see \cite{Aizenman}) that we use in this paper are described in detail in \cite{Grimm1}. The enhancement argument we utilize (see Lemma 2) is a close adaptation of that of Aizenman and Grimmett (\cite{Aizenman} or \cite{Grimm1}). We present the argument in detail here because the percolation model that we are considering does not fall into the class of models considered in \cite{Aizenman}: our model is directed and there are more than two possible states of the arrows at each site. But it does fall into the class of models treated in \cite{Holmes}.

\section{Main Results}
Let $V$ be the sub-lattice of $\mathbb{Z}^{2}$ with the sum of space and time coordinates even. At each vertex $v$ define a random variable $X(v)$ that takes values from the set of arrow configurations $\{ \ \nwarrow \ ,\ \nearrow \ , \ \nwarrow\nearrow \ , \ \odot\ \}$ with probabilities $P(X(v)=\nwarrow)$ = $P(X(v)=\nearrow)$ = $\frac{1}{2}(1-\delta)$, $P(X(v)=\nwarrow\nearrow)$ = $\delta(1-\epsilon)$, and $P(X(v)= \odot )$ = $\delta\epsilon$, where $\delta$ and $\epsilon$ are parameters with values between 0 and 1. The arrows coming out of a vertex $v$ reach to the nearby vertices in the row above to the left in case $X(v)$=``$\nwarrow$'', to the right in case $X(v)$=``$\nearrow$'' of $v$, or both in case $X(v)$=``$\nwarrow\nearrow$''. The outcome ``$\odot$'' means there are no arrows (bulk nucleation). This will define our $X$-arrow percolation model. (Note that the genealogy paths of the noisy voter model can be recovered by setting $\epsilon=1$ and for the classic voter model by setting $\delta=0$.) At 
the vertices of the same lattice place color-valued i.i.d. random variables $Y(v)$ with each color drawn uniformly out of $q$ possible colors for some finite $q\geq2$. To define the color process $Z(v)$ take the bottom half of the lattice with the top row being the row containing the origin (0,0). Direct time (row count) down, opposite of the direction of the $X$-arrows. Values for the $Z(v)$'s will be constructed using initial (time-zero) values together with the $X(v)$'s and $Y(v)$'s, as follows. At each vertex of the top (initial) row assign $Z(v)$ colors for all vertices $v$ according to some initial distribution. For each subsequent (lower) row following the initial row determine the colors of each vertex $v$ according to a rule that depends on the values of $X(v)$, $Y(v)$ and the $Z$ values of the two adjacent vertices in the row just above of $v$ called $v_{l}$ for the vertex above left of $v$ and  $v_{r}$ for the vertex above right of $v$. The rule is the following
\begin{equation*}
     Z(v)= \left \{ \begin{array}{cl}
                            Z(v_{l})    &    if\ X(v) =\ \nwarrow, \\
                            Z(v_{r})    &    if\ X(v) =\ \nearrow,\\
			    Z(v_{l})    &    if\ X(v) = \nwarrow\nearrow \ and \ Z(v_{l})=Z(v_{r}), \\
			    Y(v)   &   if\ X(v) = \nwarrow\nearrow \ and \ Z(v_{l})\neq Z(v_{r})\ or \ X(v) =\odot. \\
                  \end{array}
            \right .
\end{equation*} 
Denote by $Z_{n}$ the random variable that represents the sequence of colors $Z(v)$ of all the vertices of row $n$. $Z_{n}$ is a discrete time Markov process. Here we investigate under what conditions it is ergodic. The existence of an invariant distribution for $Z_{n}$ is a known fact that is explained in \cite{Toom} among many other sources. By ergodic we mean that the invariant distribution for $Z_{n}$ is unique and that any initial distribution converges to it.

The first question of interest is whether $X$-arrows percolate. To say it more precisely, let a path of $X$-arrows be a sequence of vertices $v_{1}$,$v_{2}$,$v_{3}$  ... such that $v_{i+1}$ lies above right or above left of $v_{i}$ with $X(v_{i})$ containing the arrow from $v_{i}$ to $v_{i+1}$. A finite path is said to terminate at $v_{n}$ when $X(v_{n})$=``$\odot$''. The $X$-arrows are said to percolate if there exists with strictly positive probability an infinite path of arrows starting at the origin or, equivalently, if with probability one there exists an infinite path of arrows starting from somewhere. 

\begin{theorem}
 $\forall$ $\delta>0$, $\exists$ $\epsilon_{c}(\delta)>0$ such that when $\epsilon>\epsilon_{c}(\delta)$, the $X$-arrows do not percolate and when $\epsilon<\epsilon_{c}(\delta)$, the $X$-arrows percolate.  
\end{theorem}

As a corollary of Theorem 1 it is shown that when the $X$-arrows do not percolate the $Z_{n}$ process is ergodic for any value of $q$. Although we have not shown ergodicity for all positive $\epsilon$, the following theorem extends the interval of $\epsilon$ for which the $Z_{n}$ process is ergodic to some positive distance below $\epsilon_{c}(\delta)$ where the $X$-arrows do percolate. Thus ergodicity is not equivalent to non-percolation.

\begin{theorem}  
 $\forall$ $\delta>0$ and $q\geq2$, $\exists$ $\epsilon_{c}'(\delta)<\epsilon_{c}(\delta)$ such that for $\epsilon>\epsilon_{c}'(\delta)$, the $Z_{n}$ process is ergodic.
\end{theorem} 
To prove Theorem 2, we will use (see Lemma 2 below) an adaptation of the percolation enhancement arguments of [1,11] that allows us to directly compare two distinct critical values for the parameter $\epsilon$: one for percolation and one for ergodicity.

A weaker result is established for all non zero $\epsilon$: 

\begin{theorem}
 $\forall \delta>0$, $\epsilon>0$, $2<q<\infty$, and any initial configuration $Z_{0}$, any sub-sequence limit distribution is color permutation invariant. 
\end{theorem}

This result does not imply convergence to a limiting distribution nor does it imply uniqueness of the invariant distribution as it is still possible to have different $k$-dimensional marginal distributions for $k\geq2$.
   
There are two special cases where the $Z_{n}$ process is ergodic for all positive $\delta$ and $\epsilon$. In the case $q = 2$, the ergodicity of $Z_{n}$ is known. This follows from the fact that for $q=2$, at ``$\nwarrow\nearrow$'' points when the colors of neighbors disagree choosing color 1 or 2 uniformly at random is equivalent to choosing the left or right neighbor uniformly at random. Therefore, the paths of arrows that are used to trace the color genealogy of each vertex are distributed as coalescing simple symmetric random walks with dying (``$\odot$'' points). Such paths are almost surely finite, and ergodicity will follow by the argument of Corollary 1 below. For $q = \infty$, instead of colored vertices the state space of $Z_{n}$ may be considered to be partitions of $\mathbb{Z}$ into infinitely many possible equivalence classes that play the role of colors. In this special case  the ergodicity is established by Theorem 4. 

\begin{theorem}
 When $q=\infty$, $\forall\delta>0$ and $\epsilon>0$, the $Z_{n}$ process is ergodic. 
\end{theorem}

The conjecture of ergodicity for $2<q<\infty$ for all $\delta>0$ and \textit{any} $\epsilon>0$ remains open.

\section{Proofs}
First we clarify a possible source of confusion: the time (row count) for paths of arrows is in the direction of the arrows which is upward, but the time  for the color evolution process $Z_{n}$ is against the arrows or downward. This is due to the color evolution rules that use arrows pointing to the colored vertices from the past in determining colors for the vertices in the present--i.e.,the arrows track the color genealogy. Throughout the rest of the paper we will choose the direction of time depending on whether we consider colors or not. 
  
\textbf{Proof of Theorem 1:}
Fix $\delta>0$. For simplicity let $x$ describe the horizontal coordinate of the lattice $V$ and $y$ describe the vertical coordinate of the lattice. Let $\Theta(\epsilon)$ be the probability of having an infinite path from the origin. The fact that $\Theta(\epsilon)$ decreases as $\epsilon$ increases follows from a coupling argument similar to the one used in the standard percolation model as we now explain. For each vertex $v$ of the lattice let $U(v)$ be a collection of i.i.d. uniform[0,1] random variables. For each $\epsilon>0$ set 
\begin{equation*}
     X_{\epsilon}(v) = \left \{ \begin{array}{cl}
                           \nwarrow     &    if\ U(v)\in[0,\frac{1-\delta}{2}), \\
                           \nearrow     &    if\ U(v)\in[\frac{1-\delta}{2},1-\delta), \\
			   \nwarrow\nearrow   &  if\ U(v)\in[1-\delta, 1-\delta\epsilon), \\
			   \odot &   if\ U(v)\in[1-\delta\epsilon,1].\\
                  \end{array}
            \right .
\end{equation*}
For a fixed $\epsilon$, the marginal distribution of $X_{\epsilon}(v)$ is the same as the distribution of the arrow-valued random variable $X(v)$ with the parameters $\delta$ and $\epsilon$. For $\epsilon_{1}<\epsilon_{2}$, $X_{\epsilon_{1}}$ has all the arrows of $X_{\epsilon_{2}}$. Hence, $\Theta(\epsilon_{1})\geq \Theta(\epsilon_{2})$. To show that $\epsilon_{c}<1$ we observe that the number of vertices at height $n$ reached from the origin by paths of arrows is dominated by the number of vertices at height $n$ of the branching process with the same distribution of offspring as the number of arrows $X(v)$ coming out of a vertex $v$. It follows that since $E(\hbox{\# of\ offspring}) = 1\cdot(1-\delta) + 2\cdot\delta(1-\epsilon) + 0\cdot\delta\epsilon < 1$ for $\epsilon>\frac{1}{2}$ then the branching process eventually almost surely dies out, which in turn implies that $\Theta(\epsilon)=0$ for such $\epsilon$.

\begin{figure}[htps]
\begin{center}
	\includegraphics[width=4in, height=2in]{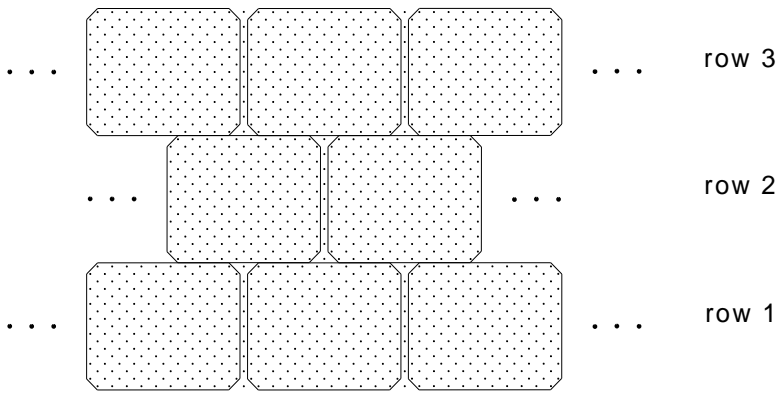}
	\caption{Relative positions of boxes.}
	\label{fig1}
\end{center}
\end{figure}

It remains to prove that $\epsilon_{c}>0$. The idea of the proof (see Lemma 1 below) is taken from an argument of Durrett \cite{Durr1} giving an upper bound for the critical probability of the oriented independent site percolation model. It is shown in \cite{Durr1} that for oriented independent site percolation on $\mathbb{Z}^{2}$, $p_{c} < 0.819$. We use a dynamic renormalization technique \cite{Barsky}. The idea behind it is described in \cite{Grimm1}. Define a box $B(k,n)$ to be $6k$ wide and $n$ high, where $k$ and $n$ will be chosen later. Such a box will be a renormalized site in an oriented independent site percolation model. Let the lower, upper, and side boundaries of the box be rows or columns of vertices parallel to and just inside the edges of the box. Position the boxes in rows such that the boxes of the first row have their lower boundary on the $t = 1$ line (here $t$ goes up) with a column of vertices between adjacent boxes' side boundaries. The second row of boxes will have its lower boundary 
right above the upper boundary of the first row but shifted by $3k+1$ to the side so that the middle of each box in the second row is above a column of vertices left in between boxes in the first row. The third row will be positioned on top of the second row similarly, and so on. See Fig. \ref{fig1} for a diagram.

\begin{figure}[htps]
\begin{center}
	\includegraphics[width=3in, height=2in]{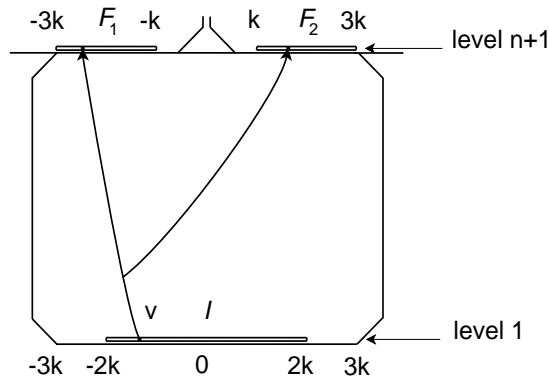}
	\caption{Box is occupied given it is checked.}
	\label{fig2}
\end{center}
\end{figure}

Since the boxes do not overlap, the events described in terms of the vertices inside of a box are independent for different boxes. For each box we define the following events. Divide the lowest row of vertices in the box into six equal parts with the length of each part equal to $k$ and place the $x=0$ coordinate in the middle of the row. For every vertex $v$ in the interval $I$ = [$-2k, 2k$] of the lowest row of vertices of the box let $A_{v}$ be the event that there exist paths from $v$ to at least one of the vertices in the interval $F_{1}$ = [$-3k, -k$] and to at least one of the vertices in the interval $F_{2}$ = [$k, 3k$] in the horizontal row of vertices right above the box (which belong to the boxes in the next row of boxes). Both paths have to be completely inside the box containing $v$. Note two things: $i$) even though the paths reach one level above the box containing $v$, their existence depends on the arrow states of vertices inside the box only; $ii$) the destination intervals $F_{1}$ and $F_{
2}$ are in corresponding intervals $I$ of the two boxes in the row of boxes above to the left and to the right of the box containing $v$. For a diagram of the described event see Fig. \ref{fig2}.

\begin{lemma}
If $\min_{v \in I} P(A_{v}) > p_{c}$, the critical value for oriented independent site percolation on $\mathbb{Z}^{2}$, then the $X$-arrows percolate. 
\end{lemma}

\textbf{Proof of Lemma 1:}
The lemma follows by comparing boxes to sites in the oriented site percolation model using the standard dynamic renormalization technique described in \cite{Barsky} and \cite{Grimm1} among other places. The idea behind the technique is to construct the infinite cluster of the origin sequentially. Consider the box containing the origin occupied if $A_{0}$ holds. Then in a predetermined order that increases with height check the boxes that are upper-right and upper-left neighbors of already occupied boxes. A box is checked if one of its neighbors below-left or below-right or both are occupied. By an inductive argument, this implies that the interval $I$ of the box being checked contains at least one vertex $v$ with a path of $X$-arrows from the origin reaching it. Choosing one such $v$ we declare the box occupied if $A_{v}$ holds. Under the assumption of Lemma 1, the cluster of occupied boxes will stochastically dominate the cluster of vertices reached from the origin in the oriented independent site 
percolation model. $\Box$ 

\begin{figure}[htps]
\begin{center}
	\includegraphics[width=2in, height=2in]{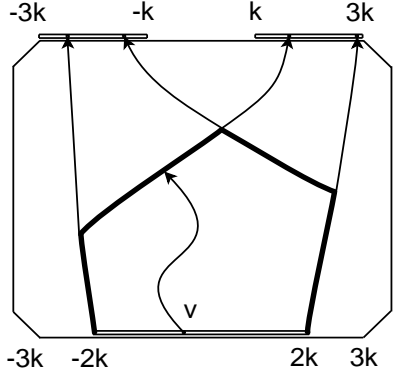}
	\caption{$A_{-2k}\cap A_{2k} \subset A_{v}$ for $-2k \leq v \leq 2k$.}
	\label{fig3}
\end{center}
\end{figure}

To apply Lemma 1 we will choose appropriate values for the height of a box $n$ and its width $6k$ that make the probability of the event $A_{v}$ arbitrarily close to one for all possible starting vertices $v$ in the interval I = [$-2k, 2k$]. First, set the parameter $\epsilon$ to 0, which sets the distribution of the outgoing arrows for each vertex to $P(\nwarrow)$ = $P(\nearrow)$ = $\frac{1}{2}(1-\delta)$, $P(\nwarrow\nearrow)$ = $\delta$, and $P(\odot)$ = $0$. Pick an arbitrarily small $\zeta>0$. Now we will determine dimensions $n$ and $6k$ of a box that will make $P(A_{-2k})$ above $1-\frac{\zeta}{2}$. Let $v$ be the vertex at the left end of the interval $I$ of the box (with horizontal coordinate $-2k$ relative to the middle of the box and vertical position $1$ relative to the bottom of the box). Note that with $\epsilon$ set to zero, at each time step we equally likely advance either to the right or to the left and sometimes both. Let $\hat{S}_{1}$ denote a path starting at this $v$ which follows the 
arrows and makes a random choice at a site with both arrows. If we consider the rightmost path $\hat{S}_{2}$ from $v$, we advance to the right with probability $P(\nearrow)$ + $P(\nwarrow\nearrow)$ = $\frac{1}{2}(1-\delta)+\delta$ = $\frac{1}{2}(1+\delta)$ and we advance to the left with probability $P(\nwarrow)$ = $\frac{1}{2}(1-\delta)$. By translating $v$ to the space-time origin we obtain coupled simple random walks $S_{1}$ and $S_{2}$ starting at the origin with $S_{1}$ symmetric and $S_{2}$ asymmetric. A corresponding translate of the event $A_{v}$ then contains the intersection event $\{ S_{1}$ $hits$ $the$ $interval$ $[-k,k]$ $at$ $time$ $n$ $without$ $going$ $below$ $-k\}$ $\bigcap$ $\{ S_{2}$ $hits$ $the$ $interval$ $[3k,5k]$ $at$ $time$ $n$ $without$ $going$ $above$ $5k\}$. Choose $k = \frac{n\delta}{4}$. The probabilities of the events defined in terms of the random walks can be easily approximated using the central limit theorem: 
\begin{equation*}
P(S_{1}(n)\in[-k,k]) = P(\frac{S_{1}}{\sqrt{n}}\in[\frac{-k}{\sqrt{n}},\frac{k}{\sqrt{n}}]) \approx P(Z\in[-\frac{\sqrt{n}\delta}{4},\frac{\sqrt{n}\delta}{4}]),
\end{equation*}
and the reflection principle:
\begin{equation*}
P(S_{1}(i) \geq -k \ for\ all \ 0 \leq i \leq n) = 1-2P(S_{1}(n)\leq-k) \approx 1-2P(Z\leq-\frac{\sqrt{n}\delta}{4}),
\end{equation*}
where $Z$ is a standard normal random variable. Similarly for $S_{2}$, whose increments $X_{i}$ have $E(X_{i})=\delta$ and $Var(X_{i})=1-\delta^{2}$, we have 
\begin{equation*}
P(S_{2}(n)\in[3k,5k]) \approx P(Z\in[-\frac{\sqrt{n}\delta}{4\sqrt{1-\delta^{2}}},\frac{\sqrt{n}\delta}{4\sqrt{1-\delta^{2}}}])
\end{equation*}
and 
\begin{equation*}
P(S_{2}(i) \leq 5k \ \forall \ 0 \leq i \leq n) \geq 1-2P(S_{2}(n) \geq 5k) \approx 1-2P(Z\geq \frac{\sqrt{n}\delta}{4\sqrt{1-\delta^{2}}}).
\end{equation*}
All these probabilities can be made greater than $1-\frac{\zeta}{8}$ by taking $n$ large enough so that the probability of their intersection is at least $1-\frac{\zeta}{2}$.
Once we have $P(A_{-2k})$, and by symmetry $P(A_{2k})$, above $1-\frac{\zeta}{2}$, since there is at least one arrow coming out of every vertex, it is easy to see from Fig. \ref{fig3} that for any $v\in[-2k,2k]$, $P(A_{v})$ $\geq$ $P(A_{-2k} \bigcap A_{2k})$ $\geq$ $1-\zeta$. (As noted by M. Damron, since also $P(\bigcap_{v=-2k}^{v=2k} A_{v})$ $\geq$ $P(A_{-2k} \bigcap A_{2k})$ $\geq$ $1-\zeta$, one can replace our dynamic by a static renormalization argument.) 

Now, as we have the size of a box fixed by choosing $n$ large enough, the probabilities of the events defined in terms of the arrows inside of a specific box are polynomials in $\epsilon$  and, therefore, continuous in $\epsilon$. If the values for $P(A_{v})$ for all the vertices $v$ in $[-2k,2k]$ exceed $1-\zeta$ for $\epsilon=0$, then they exceed $1-2\zeta$ for some small $\epsilon>0$ by continuity. This completes the proof of Theorem 1. $\Box$ 
\vspace{8mm}

One consequence of Theorem 1 is that, for fixed $\delta$ and $\epsilon>\epsilon_{c}(\delta)$, the cluster of arrows coming out of vertices of some finite interval of the initial row is almost surely finite. All the paths from those finitely many vertices terminate at ``$\odot$'' points. For the color process $Z_{n}$ this means that with the initial time long enough in the past the distribution of colors of any finite interval at time 0 for a fixed $\delta$ and $\epsilon>\epsilon_{c}(\delta)$ is determined by the distribution of a finite cluster of arrows coming out of the interval and by the i.i.d. $Y(v)$ colors at the end (``$\odot$'') points of the cluster. This proves the following corollary.

\begin{corollary}
 The $Z_{n}$ process is ergodic for any $\delta>0$ and $\epsilon>\epsilon_{c}(\delta)$. 
\end{corollary}
\vspace{8mm}

\begin{figure}[htps]
\begin{center}
	\includegraphics[width=3in, height=0.9in]{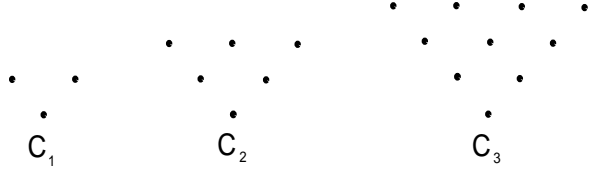}
	\caption{The sets of vertices C1, C2 and C3.}
	\label{fig4}
\end{center}
\end{figure}

To prove Theorem 2 we will use an auxiliary lemma which is an adaptation of the enhancement results of \cite{Grimm1} for standard bond percolation. The lemma is for the arrow configurations and does not use any color-valued variables. First we describe the types of enhancement we will use. On the lattice $V$ let $C_{n}$ be the deterministic subset of $V$ which consists of all the vertices up to height $n$ that can be reached from the origin by arrows taken from the set of all possible arrow configurations $\{\nwarrow,\nearrow,\nwarrow\nearrow,\ \odot\ \}$ at each vertex $v$ (see Fig. \ref{fig4}).

\begin{figure}[htps]
\begin{center}
	\includegraphics[width=4in, height=2in]{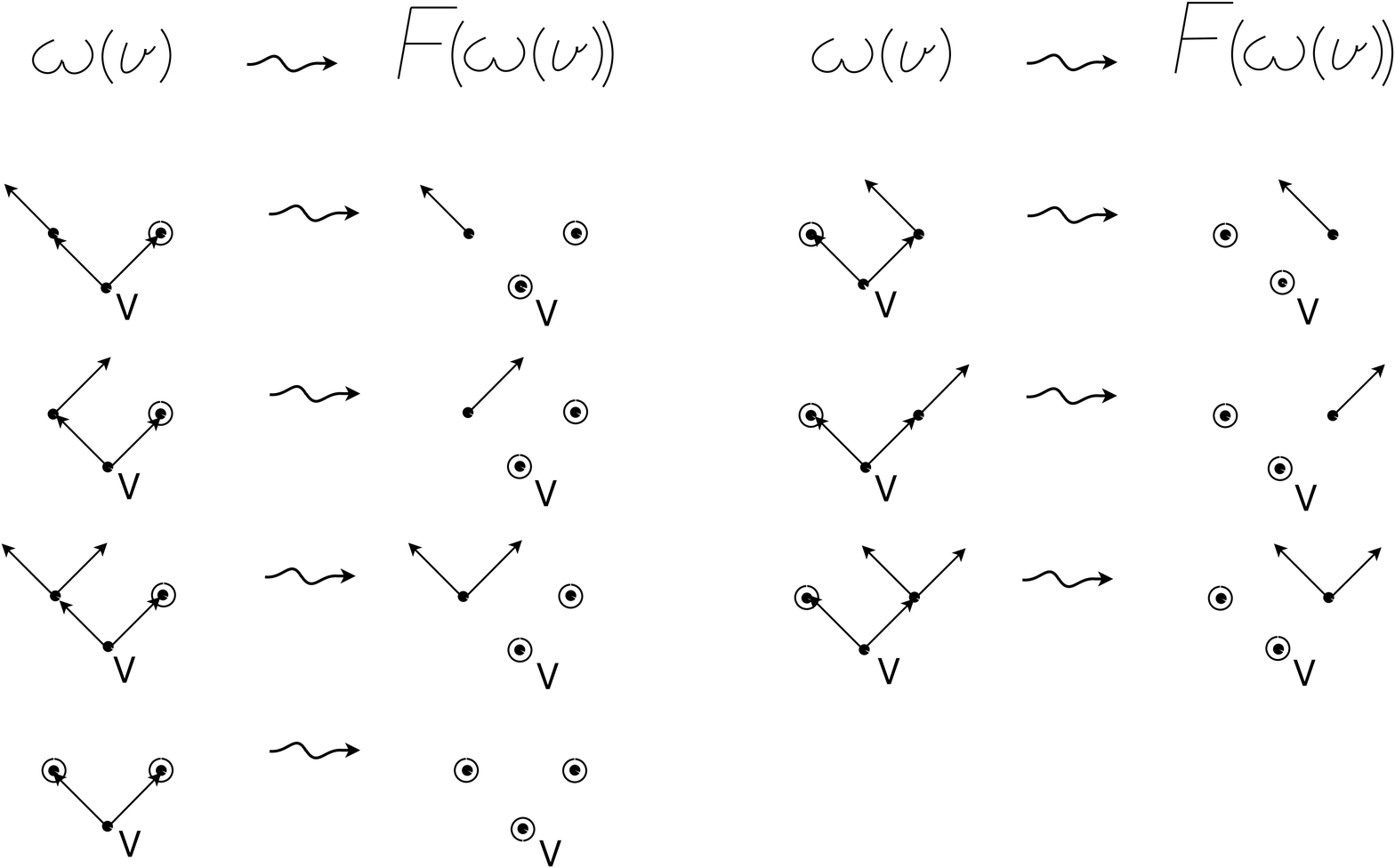}
	\caption{Configurations changed by enhancement function.}
	\label{fig5}
\end{center}
\end{figure}

Let $\Omega = \{\nwarrow,\nearrow,\nwarrow\nearrow,\ \odot\ \}^{V}$ , and $\Lambda = \{0,1\}^{V}$. For any $v\in V$ and $\omega \in \Omega$, let $F(\omega(v))$ be an enhancement function, where $\omega(v)$ is a restriction of $\omega$ to the subgraph $v\ +\ C_{1}$. $F(\omega(v))$ is defined as follows: if the arrow value at $v$ is ``$\nwarrow\nearrow$'' and at least one of the two remaining vertices in $v + C_{1}$ has value ``$\odot$'', then the arrow value at $v$ is changed to ``$\odot$'' while arrow values at other sites are unchanged; in any other case no change is made --- see Fig. \ref{fig5}.

For any configuration $\omega \in \Omega$ and $\lambda \in \Lambda$, at each vertex $v \in V$, the enhancement $\omega(v) \rightarrow F(\omega(v))$ is made if $\lambda(v) = 1$. Fix $\delta>0$. The probability measure we take on $\Omega$ is the product measure joint distribution of the previously described i.i.d. $X(v)$ random variables at the vertices $v \in V$. The probability measure we take on $\Lambda$ is the product measure joint distribution of i.i.d. Bernoulli random variables at the vertices $v$ with probability of ``$1$'' equal to $s$ and probability of ``$0$'' equal to $1-s$. The probability measure on $\Omega\times\Lambda$ denoted by $P_{\epsilon,s}$ is the product of these two measures on $\Omega$ and $\Lambda$.
Let $G^{enh}(\omega,\lambda)$ be the graph of arrows from the space $\Omega\times\Lambda$ that is obtained from the random graph $G(\omega)$ of arrows of $\Omega$ (whose marginal distribution is the same as the original arrow process distribution) by applying the enhancement function $F(\omega(v))$ at each $v$ with $\lambda(v) = 1$.
Define
\begin{equation*}
\epsilon_{c}^{enh} = \inf\{\epsilon: P_{\epsilon,s}(G^{enh}(\omega,\lambda)\ has\ an\ infinite\ path\ from\ the\ origin) = 0\} 
\end{equation*}
and
\begin{equation*}
\epsilon_{c} = \inf\{\epsilon: P_{\epsilon,s}(G(\omega)\ has\ an\ infinite\ path\ from\ the\ origin) = 0\}.
\end{equation*}
Note that $\epsilon_{c}$ is independent of $s$ and is the same as $\epsilon_{c}^{enh}$ with the parameter $s$ set to 0.
We are ready to state our auxiliary result:

\begin{lemma}
$For\ any\ s > 0,\ \epsilon_{c}^{enh}<\epsilon_{c}.$
\end{lemma}

\textbf{Proof of Lemma 2:}
In the space $\Omega\times\Lambda$ define the event $A_{n}$=$\{There$ $is$ $a$ $path$ $of$ $arrows$ $in$ $G^{enh}(\omega,\lambda)$ $from$ $0$ $to$ $the$ $top$ $of$ $C_{n}\}$. Let $\Theta_{n}(\epsilon,s) = P_{\epsilon,s}(A_{n})$. For any $v \in V$, the event $\{v$ $is$ $\omega-pivotal$ $for$ $A_{n}\}$ is the collection of all $(\omega,\lambda)$ such that for $\omega=\omega'=\omega''$ off $v$, $\omega'(v)$= ``$\nwarrow\nearrow$'' and $\omega''(v)$= ``$\odot$'', $I_{A_{n}}(\omega',\lambda)\neq I_{A_{n}}(\omega'',\lambda)$ . Here, $I_{A}$ is the indicator function of the event $A$. It should be noted from the definition above that the event $\{v$ $is$ $\omega-pivotal$ $for$ $A_{n}\}$ does not depend on the value of $\omega$ at the vertex $v$ itself. Similarly, the event $\{v$ $is$ $\lambda-pivotal$ $for$ $A_{n}\}$ is the collection of all $(\omega,\lambda)$ such that $I_{A_{n}}(\omega,\lambda') \neq I_{A_{n}}(\omega,\lambda'')$ with $\lambda=\lambda'=\lambda''$ off $v$, $\lambda'(v)$ = $0$, and $\lambda''(v)$ = 
$1$.

In the remainder of the proof we show that there exists a path from $(\epsilon',s)$ to $(\epsilon'',0)$ (with $\epsilon'<\epsilon_{c}$ and $\epsilon''>\epsilon_{c}$) along which $\Theta_{n}(\epsilon,s)$ is nondecreasing for all $n$. 
First we adapt the proof of Russo's formula from \cite{Grimm1}. Set $N=|C_{n}|$ and denote the vertices in $C_{n}$ by $v_{i}$, $1\leq i \leq N$. Suppose that for every vertex $v_{i} \in C_{n}$ the value of $\epsilon$ used to describe the distribution of $X(v_{i})$ is a distinct variable $\epsilon_{i}$. For each $v_{i}$, let $U_{i}$ be a uniform[0,1] random variable with the $U_{i}$'s independent for distinct $i$. If we assign 
\begin{equation*}
    \omega(v_{i}) = \left \{ \begin{array}{cl}
                           \nwarrow     &    if\ U_{i}\in[0,\frac{1-\delta}{2}), \\
                           \nearrow     &    if\ U_{i}\in[\frac{1-\delta}{2},1-\delta), \\
			   \nwarrow\nearrow   &  if\ U_{i}\in[1-\delta, 1-\delta\epsilon_{i}), \\
			   \odot &   if\ U_{i}\in[1-\delta\epsilon_{i},1],\\
                  \end{array}
            \right .
\end{equation*}
the new model will correspond to the model on $\Omega$ described above (when $\epsilon_{i}=\epsilon$ $\forall i$). Then
\begin{align*} 
\frac{\partial\Theta_{n}(\epsilon_{1},...,\epsilon_{N},s)}{\partial\epsilon_{i}}&=\lim_{h\rightarrow0}\frac{\Theta_{n}(\epsilon_{1},...,\epsilon_{i}+h,...,\epsilon_{N},s)-\Theta_{n}(\epsilon_{1},...,\epsilon_{i},...,\epsilon_{N},s)}{h}\\
&=\lim_{h\rightarrow0}\frac{-\delta h P_{\epsilon,s}(v_{i}\ is\ \omega-pivotal\ for A_{n})}{h} \\
&=-\delta P_{\epsilon,s}(v_{i}\ is\ \omega-pivotal\ for\ A_{n}).
\end{align*}
The second equality follows from the fact that configurations that are counted in $\Theta_{n}(\epsilon_{1},...,\epsilon_{i},...,\epsilon_{N},s)$ but not in $\Theta_{n}(\epsilon_{1},...,\epsilon_{i}+h,...,\epsilon_{N},s)$ are those where every path to the top of $C_{n}$ from the origin must use one of the two arrows coming out of $v$, and the value of $U_{i}$ is between $1-\delta(\epsilon_{i}+h)$ and $1-\delta\epsilon_{i}$. To compute the value of $\frac{\partial\Theta_{n}(\epsilon,s)}{\partial\epsilon}$ we use the chain rule and get
\begin{equation}
\frac{\partial\Theta_{n}(\epsilon,s)}{\partial\epsilon} = -\delta\sum_{i=1}^{N}P_{\epsilon,s}(v_{i}\ is\ \omega-pivotal\ for\ A_{n}).
\label{eq1}
\end{equation}
A similar argument shows that 
\begin{equation}
\frac{\partial\Theta_{n}(\epsilon,s)}{\partial s} = -\sum_{i=1}^{N}P_{\epsilon,s}(v_{i}\ is\ \lambda-pivotal\ for\ A_{n}).
\label{eq2}
\end{equation}
A vertex $v$ is ``$\omega-pivotal$ $for$ $A_{n}$'' if and only if one of the following two disjoint cases holds:

$i)$ every path from the origin to height n goes through $v$;

$ii)$ the only paths from the origin to height n that do not go through $v$ go through a vertex $u(v)$ below left or below right of $v$ such that $(a)$ $\omega(u(v))$=``$\nwarrow\nearrow$'', $(b)$ $\lambda(u(v))$=1 and $(c)$ making the value of $\omega(v)$ to be ``$\odot$'' allows for the enhancement to change the value of $\omega(u(v))$ to ``$\odot$'' killing all the paths.

\begin{figure}[htps]
\begin{center}
	\includegraphics[width=4in, height=2in]{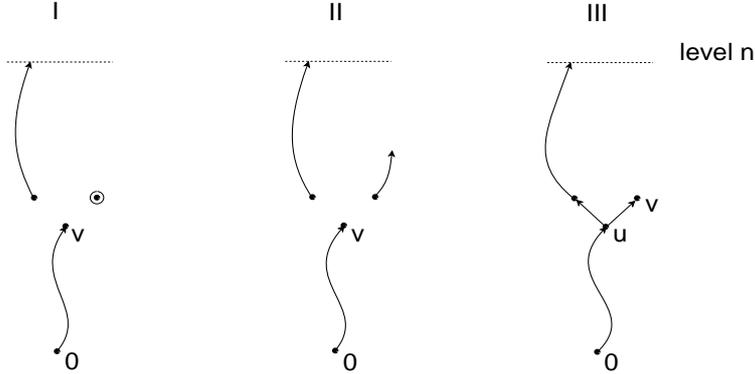}
	\caption{Examples of the configurations in the events I, II, and III.}
	\label{fig7}
\end{center}
\end{figure} 

We next decompose case  $i)$ into two disjoint events that we will denote by $I$ and $II$; case $ii)$ will be denoted $III$. Event $I$ is when the value of $\omega(v)$ fixed to be ``$\nwarrow\nearrow $'' allows for the enhancement to change it to ``$\odot$''. As seen in Fig.\ref{fig7}, event $I$ contains configurations where all possible paths from the origin to the initial level come to $v$ and continue from one of its upper-left or upper-right neighbors having the other one at state ``$\odot$''. Whether the paths exist depends, of course, on the state of $\omega(v)$, but whether $v$ is $\omega$-pivotal does not (see the definition of $\omega$-pivotal in the beginning of this proof). However, for $v$ to be $\omega$-pivotal $\lambda(v)$ must be 0, or otherwise the possible paths would not exist whether $\omega(v)$ is ``$\odot$'' (obviously) or ``$\nwarrow\nearrow$'' (would be killed by enhancement). Then by relaxing the requirement that $\lambda(v)$=``$0$'' and by fixing $\omega(v)$ to be ``$\nwarrow\
nearrow$'', we make $v$ also ``$\lambda-pivotal\ for\ A_{n}$''. 

Event $II$ is when the enhancement leaves the value of $\omega(v)$ fixed to be ``$\nwarrow\nearrow$'' unchanged (if the children of $v$ both have outgoing edges); then by switching the value of $\omega$ at the vertex above left or above right that is not involved in connecting the origin to height n (if both are involved pick the left one) to ``$\odot$'' we make $v$ to be ``$\lambda-pivotal\ for\ A_{n}$''. Therefore $\delta(1-\epsilon)\frac{1}{1-s}P_{\epsilon,s}(I)\leq P_{\epsilon,s}(v\ is\ \lambda-pivotal\ for\ A_{n})$ and $\delta(1-\epsilon)\delta\epsilon P_{\epsilon,s}(II)\leq P_{\epsilon,s}(v\ is\ \lambda-pivotal\ for\ A_{n})$.

In event $III$ (i.e. in case $ii)$), by relaxing the requirement that $\lambda(u(v))=1$ (if both below right of $v$ and below left of $v$ vertices are such as described in case $ii)$, then pick $u(v)$ to be the left one) and by fixing $\omega(v)$ to be ``$\odot$'' we make $u(v)$ to be ``$\lambda-pivotal\ for\ A_{n}$''. Therefore, $\frac{1}{s} \delta \epsilon P_{\epsilon,s}(III)\leq P_{\epsilon,s}(u(v)\ is\ \lambda-pivotal\ for\ A_{n})$. Now 
\begin{align*}
P_{\epsilon,s}&(v\ is\ \omega-pivotal\ for\ A_{n}) = P_{\epsilon,s}(I) + P_{\epsilon,s}(II) + P_{\epsilon,s}(III)\\
&\leq\frac{1-s}{\delta(1-\epsilon)} P_{\epsilon,s}(v\ is\ \lambda-pivotal\ for\ A_{n})\\
&+\frac{1}{\delta(1-\epsilon)}\frac{1}{\delta\epsilon}P_{\epsilon,s}(v\ is\ \lambda-pivotal\ for\ A_{n})\\
&+s\frac{1}{\delta\epsilon} [P_{\epsilon,s}(u_{below\ left}(v)\ is\ \lambda-pivotal\ for\ A_{n})\\
&\ \ \ \ \ \ \ \ \ \ +P_{\epsilon,s}(u_{below\ right}(v)\ is\ \lambda-pivotal\ for\ A_{n})].
\end{align*}
By summing over all the vertices in $C_{n}$ we get
\begin{align*}
\sum_{i=1}^{N}&P_{\epsilon,s}(v_{i}\ is\ \omega-pivotal\ for\ A_{n})\\
&\leq \frac{1}{\delta(1-\epsilon)}(1-s+\frac{1}{\delta\epsilon})\sum_{i=1}^{N}P_{\epsilon,s}(v_{i}\ is\ \lambda-pivotal\ for\ A_{n})\\
&+s\frac{1}{\delta\epsilon}\sum_{i=1}^{N}[P_{\epsilon,s}(u_{below\ left}(v_{i})\ is\ \lambda-pivotal\ for\ A_{n})\\
&\ \ \ \ \ \ \ \ \ \ \ \ \ \ +P_{\epsilon,s}(u_{below\ right}(v_{i})\ is\ \lambda-pivotal\ for\ A_{n})]\\
&\leq(\frac{1-s}{\delta(1-\epsilon)}+\frac{1}{\delta(1-\epsilon)}\frac{1}{\delta\epsilon}+2s\frac{1}{\delta\epsilon})\sum_{i=1}^{N}P_{\epsilon,s}(v_{i}\ is\ \lambda-pivotal\ for\ A_{n}).\\
\end{align*} 
So there is a continuous positive function $\gamma(\epsilon,s)$ such that 
\begin{align*}
\sum_{i=1}^{N}P_{\epsilon,s}(v_{i}\ is\ \omega-pivotal\ for\ A_{n})\leq\gamma(\epsilon,s)\sum_{i=1}^{N}P_{\epsilon,s}(v_{i}\ is\ \lambda-pivotal\ for\ A_{n})\
\end{align*}
and by using the Russo-like formulas \eqref{eq1} and \eqref{eq2} from above we get $\frac{\partial\Theta_{n}}{\partial\epsilon}(\epsilon,s) \geq \delta\gamma(\epsilon,s)\frac{\partial\Theta_{n}}{\partial s}(\epsilon,s)$.
Since $\gamma(\epsilon,s)$ is positive and continuous, for any $\alpha\geq 0$ we can find $M$ and $\phi \in (0,\frac{\pi}{2})$ such that $M\geq\delta\gamma(\epsilon,s)$ on $[\alpha,1-\alpha]^{2}$ and $\tan\phi = M$. The directional derivative of $\Theta_{n}(\epsilon,s)$ thus satisfies
\begin{align}
\nabla\Theta_{n}\cdot(\cos\phi,-\sin\phi) 
&= \frac{\partial\Theta_{n}}{\partial\epsilon}\cos\phi-\frac{\partial\Theta_{n}}{\partial s}\sin\phi \notag \\
&\geq\delta\gamma\frac{\partial\Theta_{n}}{\partial s}\cos\phi-\frac{\partial\Theta_{n}}{\partial s}\sin\phi \notag \\
&=-\cos\phi\frac{\partial\Theta_{n}}{\partial s}(\tan\phi-\delta\gamma) \notag \\
&\geq0 
\label{eq3}
\end{align} 

\begin{figure}[htps]
\begin{center}
	\includegraphics[width=2in, height=2in]{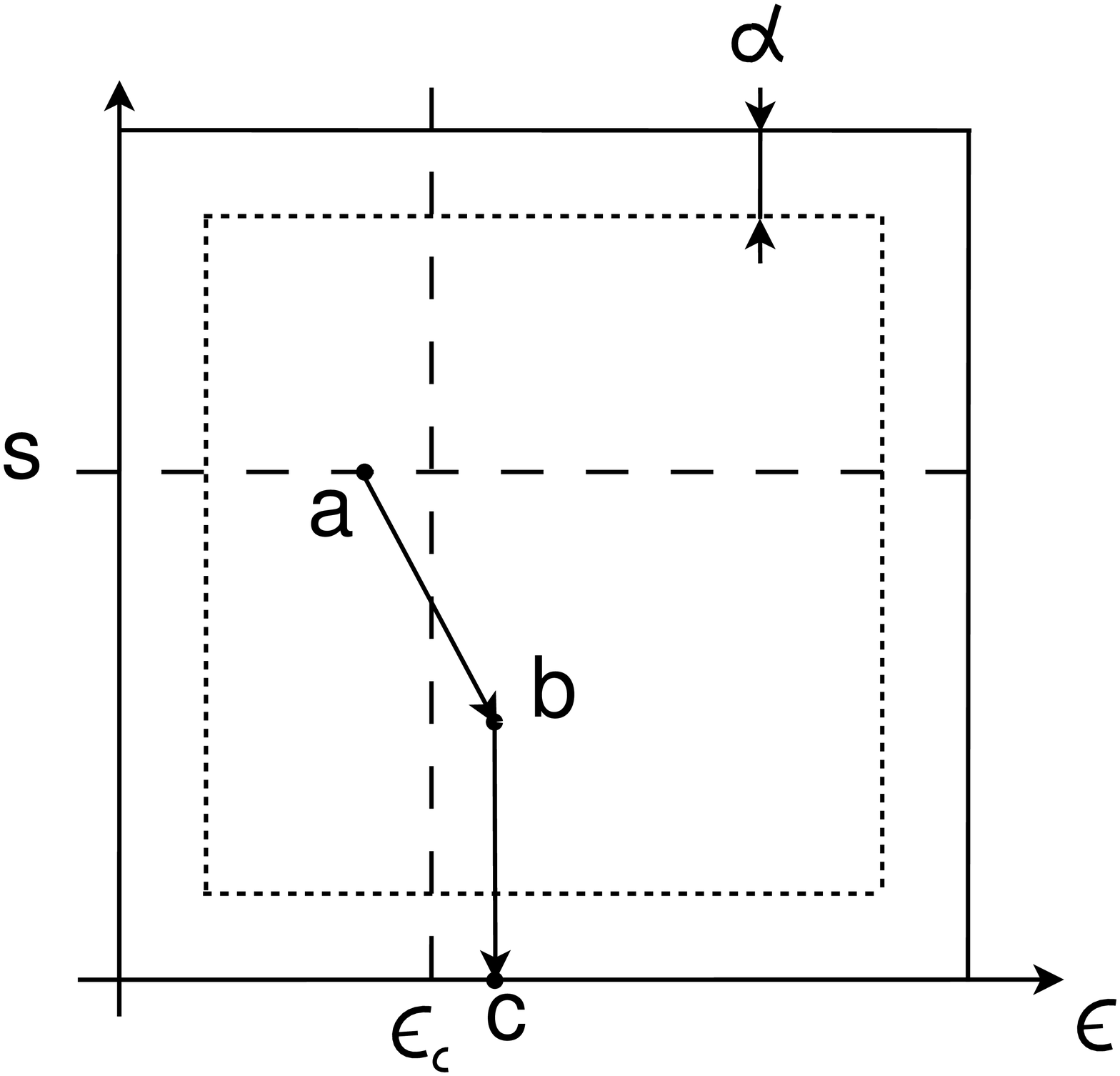}
	\caption{}
	\label{fig8}
\end{center}
\end{figure}

Set $\alpha\leq min(\frac{1}{2}\epsilon_{c},\frac{1}{2}s)$ and mark points $a$, $b$, and $c$ (see Fig. \ref{fig8}) such that $a$ and $b$ are inside of $[\alpha,1-\alpha]^{2}$, $a$ is at height $s$ but left of $\epsilon_{c}$, $b$ is such that the vector from $a$ in the direction of $(\cos\phi,-\sin\phi)$ of appropriately chosen length crosses the line $\epsilon=\epsilon_{c}$ and lands at $b$, and $c$ has the same $\epsilon$-coordinate as $b$ and has $s$-coordinate equal to $0$.
From \eqref{eq3} we conclude that $\Theta_{n}(a)\leq\Theta_{n}(b)\leq\Theta_{n}(c)$. Taking $n\rightarrow\infty$ and using that $\gamma(\epsilon,s)$ does not depend on $n$, $\Theta(\epsilon,s)$:=$P_{\epsilon,s}(There$ $is$ $an$ $infinite$ $path$ $from$ $the$ $origin)$ satisfies $\Theta(a)\leq\Theta(b)\leq\Theta(c)$. Since $\epsilon > \epsilon_{c}$ at $c$, we have $\Theta(c)=0$ and therefore $\Theta(a)=0$. This completes the proof of Lemma 2.$\Box$
\vspace{8mm}

\textbf{Proof of Theorem 2:}
We intend to couple the processes $Z_{n}$ with two different initial color sequences $A_{0}$ and $B_{0}$ in the same product space of the i.i.d. $X(v)$ and $Y(v)$ random variables to show that for any finite number of colors $q \geq 2$ and for $\epsilon$ in some interval $(\epsilon_{c}',1)$ with $\epsilon_{c}'<\epsilon_{c}$ (so that percolation of arrows occurs for $\epsilon \in (\epsilon_{c}',\epsilon_{c})$), the $Z_{n}$ process is ergodic. We note that the value of $\epsilon_{c}'$ obtained from this proof depends on the value of $q$.

Let $A$ and $B$ refer to the two processes of $Z_{n}$ corresponding respectively to the initial color sequences $A_{0}$ and $B_{0}$ (with $A(v)$ and $B(v)$ referring to the random colors of a single vertex $v$). Let $C$ be the $\{0,1\}^{\mathbb{Z}}$-valued process with $C(v) = 0$ if $A(v) = B(v)$ and $C(v) = 1$ if $A(v) \neq B(v)$. The goal is to show that almost surely $C$ eventually turns to all zeros (in any finite spatial interval). To achieve this goal we take advantage of the fact that the state space $\{0,1\}^{\mathbb{Z}}$ is partially ordered and use monotonicity present in the system (for more information see \cite{Liggett}). Let $v_{l}$ and $v_{r}$ be the vertices above left and above right of $v$. We can bound the one step transition probabilities of the $C$ process as follows. 
\begin{align*}
     P(C(v)=1|C(v_{l}),&C(v_{r}))\\
				&\leq \left \{ \begin{array}{cl}
				      0 & if\ C(v_{l}) = C(v_{r}) = 0, \\
				      P(\nearrow) + P(\nwarrow\nearrow)\frac{q-1}{q} & if\ C(v_{l}) = 0,C(v_{r}) = 1,\\
				      P(\nwarrow) + P(\nwarrow\nearrow)\frac{q-1}{q} & if\ C(v_{l}) = 1,C(v_{r}) = 0.\\
				      1-P(\odot) & if\ C(v_{l}) = C(v_{r}) = 1.\\
					      \end{array}
				      \right .
\end{align*} 
This is so because when $C(v_{l}) = C(v_{r}) = 0$, then  $A(v_{l}) = B(v_{l})$, $A(v_{r}) = B(v_{r})$, and $A(v)$ must equal $B(v)$. The case $C(v_{l}) = 0$ and $C(v_{r}) = 1$ means $A(v_{l}) = B(v_{l})$ and $A(v_{r}) \neq B(v_{r})$. To have $A(v) \neq B(v)$ we need ``$\nearrow$'' at $v$ or ``$\nwarrow\nearrow$'' at $v$ in which case we know that at least one of $A(v_{l}) \neq A(v_{r})$ or $B(v_{l}) \neq B(v_{r})$ happens (because otherwise we would have $C(v_{l}) = C(v_{r})$) and the probability that the newly chosen color is different from the color of the other process is at most $\frac{q-1}{q}$ (it is 0 if both $A(v_{l}) \neq A(v_{r})$ and $B(v_{l}) \neq B(v_{r}))$. The case $C(v_{l}) = 1$ and $C(v_{r}) = 0$ is similar. For the case $C(v_{l}) = 1$ and $C(v_{r}) = 1$, even though it is still possible that ${C(v) = 0}$ when $X(v)$ = ``$\nwarrow\nearrow$'', $A(v_{l}) \neq A(v_{r})$, $B(v_{l}) \neq B(v_{r})$, and the value of the newly chosen color $Y(v)$ goes to both $A(v)$ and to $B(v)$, we bound the 
transition probability by $1 - P(\odot)$ for simplicity. 
We keep the initial configuration of the $C$ process for a newly defined $\{0,1\}^{\mathbb{Z}}$-valued $C'$ process to which we assign the transition probabilities 
\begin{align}
\label{eq5}
     P(C'(v)=1|C'(v_{l}),&C'(v_{r})) \nonumber \\
				&= \left \{ \begin{array}{cl}
				      0 & if\ C'(v_{l}) = C'(v_{r}) = 0, \\
				      P(\nearrow) + P(\nwarrow\nearrow)\frac{q-1}{q} & if\ C'(v_{l}) = 0,C'(v_{r}) = 1,\\
				      P(\nwarrow) + P(\nwarrow\nearrow)\frac{q-1}{q} & if\ C'(v_{l}) = 1,C'(v_{r}) = 0, \\
				      1-P(\odot) & if\ C'(v_{l}) = C'(v_{r}) = 1. \\ 
					      \end{array}
				      \right .
\end{align} 
This modification of the transition rule from $C$ to $C'$ can be obtained with an appropriate coupling in which the only change is that some vertices switch from 0 to 1. An appropriate coupling of $C$ and $C'$ is the following. The probability space for the coupled process will be the product space of the i.i.d. arrow-valued random variables $X(v)$ and i.i.d. uniform[0,1] random variabls $U(v)$ at each vertex $v \in V$. For the vertices $v$ in the initial row, $A(v)$ and $B(v)$ are assigned the values at $v$ of the color sequences $A_{0}$ and $B_{0}$ respectively. The $\{0,1\}$-valued random variables $C(v)$ and $C'(v)$ take the value 0 if $A(v) = B(v)$ and value 1 if $A(v) \neq B(v)$. For the subsequent rows, the processes $A$ and $B$ evolve according to the rules of the $Z_{n}$ process with the values for the new random color $Y(v)$ chosen according to the value of $U(v)$ by splitting the range of values $[0,1]$ into $q$ equal sub-interval and assigning the $q$ colors to each of the sub-interval in some 
predetermined order. $C(v)=0$ if $A(v) = B(v)$ and $C(v)=1$ if $A(v) \neq B(v)$ as above. For the $C'$ process we want to define the transition rules in such a way that the statement ``if $C(v)=1$, then $C'(v)=1$'' always holds and such that the marginal transition probabilities of the $C'$ process are as in $(\ref{eq5})$.

For $v$ in the initial row, the statement ``if $C(v)=1$, then $C'(v)=1$'' holds since $C(v)=C'(v)$. At each vertex $v$ after time 0 look at the values of all the processes involved at the vertices $v_{l}$ and $v_{r}$ in the previous row and at the values of $X(v)$ and $U(v)$. Suppose that in the row above which contains $v_{l}$ and $v_{r}$ it is true that whenever $C(v)=1$, then $C'(v)=1$. If $X(v)=``\odot$'', assign $C'(v) = C(v) = 0$. If $X(v) = ``\nwarrow $'' or $X(v) = ``\nearrow $'', then the value for $C'(v)$ is chosen to be $C'(v_{l})$ or $C'(v_{r})$ respectively, and the statement ``if $C(v)=1$, then $C'(v)=1$'' carries over from the previous row. In case $X(v)=``\nwarrow\nearrow$'' and $C'(v_{l})$ and $C'(v_{r})$ agree in value, $C'(v)$ takes that value. If the values of $C'(v_{l})$ and $C'(v_{r})$ disagree, the rules for choosing the new value for $C'(v)$ will depend on the values of $A(v_{l})$, $A(v_{r})$, $B(v_{l})$, and $B(v_{r})$. Unless it is a case where a special rule is required, choose the 
new value for $C'(v)$ in the similar way as the value for $Y(v)$ is chosen with value 0 assigned to the first sub-interval $[0,\frac{1}{q}]$ and value 1 assigned to the rest. If $C(v_{l}) = C(v_{r}) = 0$, then $C(v)=0$ and the statement holds. If $C(v_{l}) = C'(v_{l}) = 0$ and $C(v_{r}) = C'(v_{r}) = 1$, consider the values of $A(v_{l})$, $A(v_{r})$, $B(v_{l})$, and $B(v_{r})$. If $A(v_{l}) \neq A(v_{r})$ and $B(v_{l}) \neq B(v_{r})$, then $C(v)=0$ and the statement holds. If $A(v_{l}) \neq A(v_{r})$ and $B(v_{l}) = B(v_{r})$, the special rule is required, as follows. Choose a value for $Y(v)$ and a value for $C'(v)$ based on the outcome of the uniform[0,1] random variable $U(v)$ the following way: split the interval $[0,1]$ into $q$ equal sub-intervals and assign the $q$ colors to them in the same predetermined order as before for the $Y(v)$ value, and for the $C'(v)$ value assign 0 to the sub-interval corresponding to the color $B(v_{r})$ and 1 to the other sub-intervals. The special rule is consistent 
with the rules for choosing values at $v$ for the processes $A$, $B$, and $C$, the transition probabilities for $C'(v)$ correspond to \eqref{eq5}, and the statement ``if $C(v)=1$, then $C'(v)=1$'' holds. If $A(v_{l}) = A(v_{r})$ and $B(v_{l}) \neq B(v_{r})$ the special rule is the same except the value 0 for the $C'(v)$ is assigned to the sub-interval that holds the color of $A(v_{r})$. The case when $C(v_{l}) = C'(v_{l}) = 1$ and $C(v_{r}) = C'(v_{r}) = 0$ is similar, and in all other cases the statement holds trivially. We can say that the $C'$ process dominates the $C$ process in the sense that the statement ``if $C(v)=1$, then $C'(v)=1$'' always holds.

Therefore, if with any initial configuration the new process $C'(v)$ turns to all zeros (or ``dies'') eventually almost surely (in any finite spatial interval), then in the original color and arrow process $Z_{n}$ the colors of all the vertices (in any finite spatial interval) will eventually become the same for all initial color configurations, which implies the uniqueness of the invariant distribution and convergence to it.

To show that the $C'$ process eventually turns to all zeros, we use the auxiliary enhancement result of Lemma 2. To relate the random graph $G$ of $X(v)$ arrows before and after the enhancement to a $\{0,1\}^{V}$-valued processes, we think of all the vertices of $G$ (or $G^{enh}$) that have an infinite path coming out of them as having value 1 and all the vertices whose (directed) cluster of arrows is finite as having value 0. With this in mind, for illustration purpose we first relate the process of $X$-arrows before the enhancement to a new $\{0,1\}^{V}$-valued process $C''$ defined by assigning $C''(v)$ = 1 to all vertices at the initial time, by directing time down (opposite of the arrow direction), and by transferring the value of 1 down by any arrow pointing to a site with value 1. The transition probabilities of $C''$ are the following:
\begin{align*}
     P(C''(v)=1|C''(v_{l}),&C''(v_{r}))\\
				&= \left \{ \begin{array}{cl}
				      0 & if\ C''(v_{l}) = C''(v_{r}) = 0, \\
				      P(\nearrow) + P(\nwarrow\nearrow) & if\ C''(v_{l}) = 0,C''(v_{r}) = 1,\\
				      P(\nwarrow) + P(\nwarrow\nearrow) & if\ C''(v_{l}) = 1,C''(v_{r}) = 0,\\
				      1-P(\odot) & if\ C''(v_{l}) = C''(v_{r}) = 1.\\
					      \end{array}
				      \right .
\end{align*} 
The $C''$ process can be defined independently on $\{0,1\}^{V}$ or it can be coupled with the process of $X$-arrows as described above with the same transition rates. An almost sure limiting configuration of all zeros of the $C''$ process is equivalent to non-percolation of arrows. To relate the enhanced $X$-arrows process with a $\{0,1\}^{V}$-valued process we consider the $C''$ process realized on the enhancement space $\Omega\times\Lambda$ introduced before Lemma 2. The configuration of $C''$ depends only on the graph of arrows $G(\omega)$. Now, let $C^{*}$ denote the $\{0,1\}^{V}$-valued enhanced percolation of arrows process that is constructed from $G^{enh}(\omega,\lambda)$ the same way $C''$ is constructed from  $G(\omega)$: $C^{*}(v)=1$ if $v$ is connected by arrows to time 0 in $G^{enh}(\omega,\lambda)$, and $C^{*}(v)=0$ otherwise. The enhancement function described in the arrow percolation context is activated with probability $s$ at a vertex $v$ only when we have at $v$ a double-arrow and at least one of the 
vertices above right or above left of $v$ has ``$\odot$'' value. For the $C^{*}$ realization to activate the enhancement it is necessary that at least one of the vertices above right or above left (the one with $X$-value = $\odot$) has value $0$. But that is not sufficient since it is possible that both $v_{l}$ and $v_{r}$ have arrows coming out of them but with all the paths coming out of, for example, $v_{l}$ terminating at ``$\odot$'' points. Therefore, the transition probabilities for $C^{*}$ process will satisfy the following:
\begin{align*}
     P(C^{*}(v)=1|C^{*}&(v_{l}),C^{*}(v_{r}))\\
					     &\geq \left \{ \begin{array}{cl}
					      0 & if\ C^{*}(v_{l}) = C^{*}(v_{r}) = 0, \\
					      P(\nearrow) + P(\nwarrow\nearrow)(1-s) & if\ C^{*}(v_{l}) = 0,C^{*}(v_{r}) = 1,\\
					      P(\nwarrow) + P(\nwarrow\nearrow)(1-s) & if\ C^{*}(v_{l}) = 1,C^{*}(v_{r}) = 0,\\
					      1-P(\odot) & if\ C^{*}(v_{l}) = C^{*}(v_{r}) = 1.\\
							    \end{array}
						  \right .
\end{align*}
Set $s = \frac{1}{q}>0$. For that choice of $s$ with the standard coupling we can dominate the process $C'$ by the process $C^{*}$ in the sense that ``if $C'(v)=1$, then $C^{*}(v)=1$''. Lemma 2 says that $\epsilon_{c}^{enh}$=$\epsilon_{c}' < \epsilon_{c}$. This means that for the interval of $\epsilon_{c}'<\epsilon \leq\ 1$ the enhanced cluster of the origin is almost surely finite, and all the values come from ``$\odot$'' points which give the value $0$. Therefore, the unique invariant distribution for such $\epsilon$ for the processes $C^{*}$, $C'$, and $C$ is the delta measure at all zeros.$\Box$
\vspace{8mm}

\textbf{Remark:}
The dominating $C'$ process considered in the proof of Theorem 2 is a simpler $\{0,1\}^{\mathbb{Z}}$-valued process for which many techniques are developed in discrete and continuous time (\cite{Toom,Griffeath,Liggett} among others). For example, Theorem 4.1 of Chapter I or the proof of Theorem 3.32 of Chapter III of \cite{Liggett} provide for an easily derived lower bound such that for $\epsilon > \frac{q-2}{2q-2}$ the continuous time version of the $C'$ process a.s. converges to the all zero configuration from any initial state. The proof of Theorem 3.32 easily extends to the discrete time version of $C'$ providing the same lower bound. However, we chose to use the enhancement argument of Lemma 2 because it allows us to compare the actual critical values $\epsilon_{c}$ for percolation and $\epsilon_{c}'$ for ergodicity.    
\vspace{8mm}

\begin{figure}[htps]
\begin{center}
	\includegraphics[width=4in, height=1.2in]{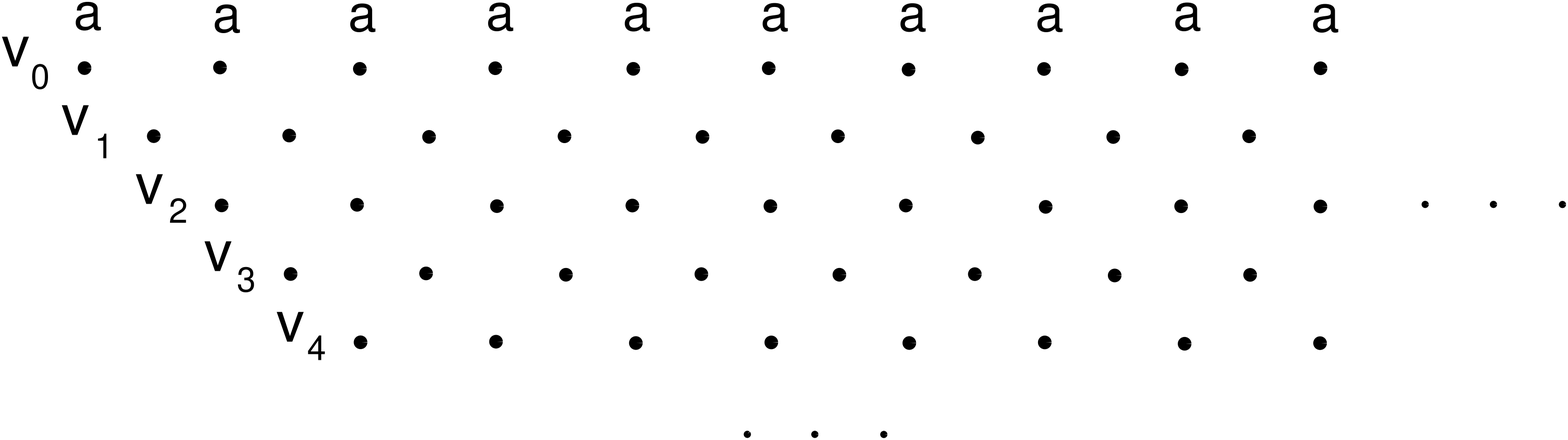}
	\caption{Set $V'$.}
	\label{fig11}
\end{center}
\end{figure}

\textbf{Proof of Theorem 3:}
Let $V' \subset V$ be the set of vertices from which only the horizontal half-line  to the right of and including the origin can be reached by arrows taken from the set of all possible arrow configurations $\{\nwarrow,\nearrow,\nwarrow\nearrow,\odot \}$ at each vertex  $v \in V$. See Fig \ref{fig11}.

On the space $\Omega=\{\nwarrow,\nearrow,\nwarrow\nearrow,\odot \}^{V}\times\{a, b, c, ... \}^{V}$, where $a, b,...$ are $q$ colors, let $P$ be the product distribution of the i.i.d. arrow-valued $X(v)$ and the color-valued uniform $Y(v)$ random variables for $v \in V$. For each $\omega \in \Omega$ and for each $v \in V'$ define new arrow-valued random variables $W(v)(\omega)$ and color-valued random variables $Z(v)(\omega)$ as follows: for $v$ in the top row of $V'$ fix $W(v)(\omega)=$ ``$\odot$'' and assign some configuration of colors $Z(v)$ to the top row. For the vertices in each consecutive row below the top one define
\begin{equation*}
     W(v)(\omega)= \left \{ \begin{array}{cl}
                           \ \nwarrow     &    if\ X(v)(\omega) =\ \nwarrow, \\
                           \ \nearrow     &    if\ X(v)(\omega) =\  \nearrow, \\
			   \nwarrow\nearrow   &  if\ X(v)(\omega) = \nwarrow\nearrow\ and\  Z(v_{l})(\omega)=Z(v_{r})(\omega), \\
			   \ \odot &   otherwise, \\
                  \end{array}
            \right .
\end{equation*} 
(here $v_{l}$ means the vertex above left of $v$ and $v_{r}$ the vertex above right of $v$) and
\begin{equation*}
     Z(v)(\omega)= \left \{ \begin{array}{cl}
                            Z(v_{l})(\omega)    &    if\ W(v)(\omega) =\ \nwarrow, \\
                            Z(v_{r})(\omega)    &    if\ W(v)(\omega) =\ \nearrow, \\
			    Z(v_{l})(\omega)=Z(v_{r})(\omega)  &  if\ W(v)(\omega) = \nwarrow\nearrow, \\
			    Y(v)(\omega)   &   if\ W(v)(\omega) =\ \odot. \\
                  \end{array}
            \right .
\end{equation*} 

Thus defined, $Z(v)$ is the same as the color of the vertex $v$ of the $Z_{n}$ process. We want to estimate how the probability that the color $Z(v)$ comes from the top row changes. Let $G$ be the collection of vertices $v \in V'$ that are connected by paths of $W$-arrows to the top row. 
Let $v_{0},v_{1},...$ be the leftmost vertices of $V'$ as in Fig. \ref{fig11} and consider the sum 
\begin{align*}
\sum_{i=0}^{\infty} P(v_{i} \in G).
\end{align*}
Denote the upper right vertex of $v_{i+1}$ by $u_{i}$ and the event $\{At$ $least$ $one$ $of$ $v_{i}$ $or$ $u_{i}$ $is$ $in$ $G$ $and$ $the$ $other$ $one$ $has$ $the$ $same$ $color\}$ by $E_{i}$. Then 
\begin{align}
\label{eq4}
P(v_{i+1} \in G) =& P(X(v_{i+1}) = \nwarrow \ )P(v_{i} \in G) + P(X(v_{i+1}) = \nearrow \ )P(u_{i} \in G) \nonumber \\
&+P(X(v_{i+1}) = \nwarrow\nearrow \ )P(E_{i}), 
\end{align}
with
\begin{align*}
P(E_{i}) &= P(v_{i} \in G\ and\ u_{i} \in G \ and\ Z(u_{i}) = Z(v_{i}))\\
&+ P(v_{i} \in G\ and\ u_{i} \notin G\ and\ Z(u_{i})=Z(v_{i}))\\
&+ P(u_{i} \in G\ and\ v_{i} \notin G\ and\ Z(v_{i})=Z(u_{i}))\\
\leq P(v_{i} \in G\ and\ u_{i} &\in G ) + P(v_{i} \in G\ and\ u_{i} \notin G )\frac{1}{q}
+ P(u_{i} \in G\ and\ v_{i} \notin G )\frac{1}{q}\\
\leq P(v_{i} \in G\ and\ u_{i} &\in G ) + P(v_{i} \in G\ and\ u_{i} \notin G )\frac{1}{2}
+ P(u_{i} \in G\ and\ v_{i} \notin G )\frac{1}{2}\\
= P(v_{i} &\in G )\frac{1}{2} + P(u_{i} \in G )\frac{1}{2}
\leq \sup_{v\ in\ row\ i} P(v \in G ).
\end{align*}
Since $P(v_{i} \in G)$, $P(u_{i} \in G)$, and $P(E_{i})$ are all not greater than 
\begin{equation*}
\sup_{v\ in\ row\ i}P(v \in G ),
\end{equation*}
\begin{equation*}
P(v_{i+1}\in G)\leq(1-\delta\epsilon)\sup_{v\ in\ row\ i}P(v\in G).
\end{equation*}
Also, since the equivalent of \eqref{eq4} and the calculations above are valid
for all the vertices $v$ in row $i+1$, we have
\begin{equation*}
\sup_{v\ in\ row\ i+1}P(v \in G)\leq(1-\delta\epsilon)\sup_{v\ in\ row\ i}P(v \in G ), 
\end{equation*}
and
\begin{align*}
\sum_{i=0}^{\infty} P(v_{i} \in G) &\leq \sum_{i=0}^{\infty} \sup_{v\ in\ row\ i}P(v\in G) \\
&\leq \sup_{v\ in\ row\ 0}P(v \in G)\sum_{i=0}^{\infty}(1-\delta\epsilon)^{i} < \infty. 
\end{align*}
Therefore, for any initial coloring, almost surely $v_{i}$ is not connected by the $W$-arrows to the top of $V'$ for large enough $i$. The same is true for the $W$-arrow cluster of any deterministic finite interval of the $i^{th}$ row of vertices. We can conclude that for any subsequence limit distribution, any joint distribution of colors of a finite collection of vertices is color permutation invariant. $\Box$
\vspace{8mm}

\textbf{Proof of Theorem 4:}
We couple the $q=\infty$ process with any initial color configuration (or equivalence class partition) to the one with the initial row of a constant color. By the argument of the proof of Theorem 3, for any finite interval of vertices, if we start with the constant color configuration some sufficiently large amount of time earlier, the cluster $G$ of vertices that can be connected by the $W(v)$ arrows to the initial row will not intersect the finite interval. For $q=\infty$, we show that colors of the vertices not in the cluster $G$ obtained with the constant color initial configuration are almost surely the same for all initial color configurations. This will follow if we can show that $W(v)$ = ``$\odot$'' with the initial constant color implies $W(v)$ = ``$\odot$'' for any other initial color configuration because at such vertices a new color is chosen which is independent of the initial configuration. And this is true because, given that the initial colors are constant, $W(v)$ = ``$\odot$'' happens when 
either $X(v)$ = ``$\odot$'' and, therefore, $W(v)$ = ``$\odot$'' for any initial color configuration, or when $X(v)$ = ``$\nwarrow\nearrow$'' and the colors of $v_{l}$ and $v_{r}$ are different. If we start with the constant color configuration, for $v_{l}$ and $v_{r}$ to be of different color at least one of them, say $v_{r}$, has to be connected by the $X$-arrows to some vertex $w$ that is located later than the initial time and such that $X(w)$ = ``$\odot$'' (or both, $v_{l}$ and $v_{r}$, have to be connected to distinct such $w$). For the case $q = \infty$ (when every ``$\odot$'' vertex almost surely creates a new color that does not match any previously existing color), every ``$\nwarrow\nearrow$'' along the path of the $X(v)$ arrows from $v_{r}$ to $w$ will be changed to ``$\odot$'' (unless the double arrow creates branches that come together again before reaching $w$). But in this situation the same thing will happen along the path from $v$ to $w$ for any initial color configuration, and colors of the 
vertices above left and above right of $v$ will also be different no matter what the initial color configuration was. Therefore $W(v)$ = ``$\odot$'' for every initial color configuration. Therefore, if starting from the constant color initial configuration the distribution of equivalence class partition at time $n$ converges as $n\rightarrow\infty$, then the $q=\infty$ system is ergodic.

The fact that the limit distribution of the equivalence class partition exists as $n\rightarrow\infty$ for a constant color initial configuration follows from the argument above because it allows us to couple equivalence class partitions corresponding to different initial times. Instead of looking at the distributions of the equivalence class partitions of a finite interval of vertices at time $n$ and time $n+1$ corresponding to the constant color initial configuration, we look at the interval at time 0 and compare the distributions of its equivalence class partitions corresponding to the initial times $-n$ and $-n-1$. The equivalence class partition of the interval corresponding to the initial time $-n-1$ can be coupled to the one corresponding to the initial time $-n$ by adding a row of the $X$ arrows at time $-n$ that point to time $-n-1$ and moving the initial constant color configuration from row $-n$ to row $-n-1$. By the argument above, if $n$ is large enough, colors of the interval with high 
probability will be the same for all initial configurations at time $-n$ as they would be for the constant color initial configuration. Therefore, moving the constant color initial configuration from time $-n$ to time $-n-1$ will change colors at time $-n$ but not the colors of the interval. 
$\Box$

\section*{Acknowledgments}
The authors thank Mark Holmes for useful discussions. They also thank Michael Damron and an anonymous referee for a careful reading of and comments on the manuscript.

\end{document}